\documentclass[11pt]{article}

\usepackage[T1]{fontenc}
\usepackage[a4paper,margin=1in]{geometry}

\usepackage{amsmath,amsthm,amssymb,mathtools}
\usepackage{mathrsfs}
\usepackage{float}

\usepackage{bm}

\usepackage{enumitem}

\usepackage{hyperref}
\usepackage[nameinlink,noabbrev]{cleveref}

\usepackage{tikz}
\usepackage{tikz-cd}
\usetikzlibrary{arrows.meta,positioning,calc}

\usepackage[backend=biber,style=alphabetic,maxbibnames=99]{biblatex}
\newtheorem{theorem}{Theorem}[section]
\newtheorem{proposition}[theorem]{Proposition}

\newtheorem{lemma}[theorem]{Lemma}

\newtheorem{corollary}[theorem]{Corollary}

\theoremstyle{definition}

\theoremstyle{remark}
\newtheorem{remark}[theorem]{Remark}

\theoremstyle{plain}

  \Crefname{theorem}{Theorem}{Theorems}
  \Crefname{thm}{Theorem}{Theorems}
  \Crefname{lemma}{Lemma}{Lemmas}

  \Crefname{lem}{Lemma}{Lemmas}
  \Crefname{remark}{Remark}{Remarks}
  \Crefname{prop}{Proposition}{Propositions}
  \Crefname{proposition}{Proposition}{Propositions}
  \Crefname{problem}{Problem}{Problems}
\Crefname{notation}{Notation}{Notations}
\Crefname{shadethm}{Theorem}{Theorems}
\Crefname{claim}{Claim}{Claims}
  \Crefname{defn}{Definition}{Definitions}
  \Crefname{corollary}{Corollary}{Corollaries}
  \Crefname{shadedefn}{Definition}{Definitions}
  \Crefname{section}{Section}{Sections}
  \Crefname{figure}{Figure}{Figures}
  \Crefname{shadeprop}{Proposition}{Propositions}
  \Crefname{exercise}{Exercise}{Exercises}
    \Crefname{assumption}{Assumption}{Assumptions}
  
\newcommand{\Jac}{\mathrm{Jac}}

\title{\textbf{Rational Torsion and Ramification Constraints on Odd Superelliptic Jacobians}}

\author{Ankita Das\thanks{SRM University, AP, India. Email: \texttt{ankita\_das@srmap.edu.in}}}

\date{}

\begin{document}

\maketitle

\begin{abstract}
Let $C:y^d=f(x)$ be a smooth projective superelliptic curve over
$\mathbb{Q}$, where $d\geq 3$ is odd, $f(x)\in\mathbb{Q}[x]$ is
separable of degree $n\geq 3$, and $\gcd(d,n)=1$. Let
$J=\operatorname{Jac}(C)$ be its Jacobian. We study the rational
torsion subgroup $J(\mathbb{Q})_{\mathrm{tors}}$ through the
interaction of ramification geometry, Galois actions, reduction at
primes of good reduction, and the geometry of divisors. First, we
describe the subgroup $R_d\subseteq J[d]$ generated by the
ramification-point differences and determine its Galois-invariant
subgroup explicitly in terms of the sizes of the Galois orbits of
the roots of $f$. This yields a precise description of the rational
$d$-torsion classes arising from ramification, without assuming
that they account for all of $J(\mathbb{Q})[d]$. Second, we establish
criteria based on the orders of reduced Jacobians that can prove
the vanishing of the full rational torsion subgroup. Finally, we
derive geometric constraints on divisor classes represented by
effective divisors of degree smaller than the genus, including a
restriction arising when distinct effective divisors represent the
same class. Explicit examples illustrate how the Galois-theoretic
and reduction-theoretic methods can be combined to study rational
torsion on odd superelliptic Jacobians.
\end{abstract}

\tableofcontents

\section{Introduction}
\label{sec-1}

Let
\[
    C:\quad y^d=f(x)
\]
be a smooth projective superelliptic curve over $\mathbb{Q}$, where
$d\geq 3$ is odd, $f(x)\in\mathbb{Q}[x]$ is separable of degree
$n\geq 3$, and $\gcd(d,n)=1$. Let $J=\operatorname{Jac}(C)$ be its
Jacobian, and write $J(\mathbb{Q})_{\mathrm{tors}}$ for the subgroup
of rational torsion points. The study of rational torsion on
Jacobians is a central problem in arithmetic geometry: it asks
which divisor classes of finite order are defined over the ground
field, and how the geometry and arithmetic of the underlying curve
constrain their possible orders and structure.

For elliptic curves, the Nagell--Lutz theorem gives an explicit
integrality condition for rational torsion, and Mazur's theorem
classifies the possible rational torsion subgroups over
$\mathbb{Q}$ \cite{Mazur1977}. In higher genus, there is no comparable
general classification for Jacobians of arbitrary curves. One
therefore seeks effective structural information and criteria that
can be applied to specific families or individual curves. The
geometry of algebraic curves, including special divisors and
Weierstrass points, provides important tools in this setting
\cite{ACGH,Coppens1990}. The arithmetic study of Jacobians also draws
on the theory of N\'eron models and reduction, which provides a
framework for relating torsion over number fields to torsion over
finite residue fields \cite{BoschLutkebohmertRaynaud,SilvermanAdvanced}.

For superelliptic curves, the covering structure provides a natural
starting point. The map
\[
    \pi:C\longrightarrow\mathbb{P}^{1},
    \qquad (x,y)\longmapsto x
\]
is a cyclic covering of degree $d$, and its ramification points
yield distinguished divisor classes on the Jacobian. The structure
of torsion subgroups associated with superelliptic Jacobians has
been studied from several perspectives; see, for example,
Jedrzejak \cite{Jed14} and Wawrów \cite{Wawrow2021}. The cyclic
symmetries of curves and their induced actions on Jacobians also
fit into the broader study of group actions and decompositions of
Jacobians \cite{KaniRosen1989}.

Related Diophantine questions involving superelliptic curves have
also been investigated. In particular, our earlier work on a
generalized fruit Diophantine equation connects a Diophantine
problem with the study of superelliptic curves \cite{BanerjeeChakrabortyDas}. This
provides part of the broader motivation for examining how
arithmetic conditions on a curve interact with its geometric
structure. The present paper focuses specifically on rational
torsion on Jacobians, using ramification classes, Galois actions,
reduction at primes of good reduction, and restrictions arising
from low-degree divisors. For further background on the arithmetic
of curves and their Jacobians, we refer to Hindry and Silverman
\cite{HindrySilverman2000}.

The aim of this paper is to combine three complementary approaches
to rational torsion: the Galois action on ramification classes,
reduction at primes of good reduction, and the geometry of
low-degree divisors. The first approach identifies a distinguished
subgroup of the geometric $d$-torsion and computes its rational
Galois invariants. The second gives criteria that can force the
entire rational torsion subgroup to vanish, without assuming that
all rational torsion comes from ramification. The third describes
geometric restrictions on divisor classes that admit effective
representatives of degree smaller than the genus. These approaches
answer different questions, and keeping them separate is essential:
the ramification calculation alone does not classify all of
$J(\mathbb{Q})[d]$.

\medskip
\noindent
\textbf{Ramification classes and Galois invariants.}
Let $\alpha_1,\ldots,\alpha_n$ be the roots of $f$ in
$\overline{\mathbb Q}$, and let
\[
R_i=(\alpha_i,0),\qquad 1\leq i\leq n,
\]
be the corresponding geometric ramification points. Since
$\gcd(d,n)=1$, the curve has a unique point $\infty$ above infinity.
The divisor relations
\[
\operatorname{div}(x-\alpha_i)=dR_i-d\infty,
\qquad
\operatorname{div}(y)=\sum_{i=1}^nR_i-n\infty
\]
show that the classes
\[
r_i=[R_i-\infty]\in J[d](\overline{\mathbb Q})
\]
are $d$-torsion and satisfy $\sum_i r_i=0$. We define the
ramification-supported subgroup
\[
R_d=\langle r_1,\ldots,r_n\rangle
\subseteq J[d](\overline{\mathbb Q}).
\]
Although $R_d$ is stable under the absolute Galois group
$G_{\mathbb Q}=\operatorname{Gal}(\overline{\mathbb Q}/\mathbb Q)$,
its generators need not be rational. Consequently,
\[
R_d^{G_{\mathbb Q}}\subseteq J(\mathbb Q)[d],
\]
and equality is not automatic.

The first structural result identifies the abstract group generated
by these classes and gives a concrete description of its Galois
invariants. The following statements are proved in Section \ref{sec-4}.

\medskip
\noindent
\begin{lemma}
\label{lem-4.1}
[Structure of the ramification subgroup]
The homomorphism
\[
\Phi:(\mathbb Z/d\mathbb Z)^n\longrightarrow R_d,\qquad
(a_1,\ldots,a_n)\longmapsto\sum_{i=1}^n a_i r_i
\]
is surjective and has kernel equal to the diagonal subgroup
\[
\Delta=\{(t,\ldots,t):t\in\mathbb Z/d\mathbb Z\}.
\]
In particular,
\[
R_d\cong(\mathbb Z/d\mathbb Z)^n/\Delta
\cong(\mathbb Z/d\mathbb Z)^{n-1}.
\]
\end{lemma}

\medskip
\noindent
\begin{theorem}
 \label{thm-4.2} [Galois invariants of the ramification subgroup]
Suppose that the roots of $f$ decompose into $G_{\mathbb Q}$-orbits
\[
\{\alpha_1,\ldots,\alpha_n\}
=O_1\sqcup\cdots\sqcup O_s,
\qquad m_j=|O_j|.
\]
Then there is an isomorphism
\[
R_d^{G_{\mathbb Q}}\cong
\left\{(b_1,\ldots,b_s)\in(\mathbb Z/d\mathbb Z)^s:
\sum_{j=1}^s m_jb_j=0\text{ in }\mathbb Z/d\mathbb Z\right\}.
\]
Under this isomorphism, $(b_1,\ldots,b_s)$ corresponds to
\[
\sum_{j=1}^s b_j
\sum_{\alpha_i\in O_j}(R_i-\infty).
\]
\end{theorem}

This result reduces the computation of the rational classes generated
by ramification differences to a single congruence involving the
sizes of the Galois orbits. In particular, it shows precisely how the
arithmetic symmetry of the branch locus affects this distinguished
subgroup.

\medskip
\noindent
\begin{corollary}
    \label{cor-4.4}[The transitive case]
If $G_{\mathbb Q}$ acts transitively on the roots of $f$, then
\[
R_d^{G_{\mathbb Q}}\cong
\{b\in\mathbb Z/d\mathbb Z:nb=0\}.
\]
In particular, if $\gcd(d,n)=1$, then
\[
R_d^{G_{\mathbb Q}}=0.
\]
This is a statement about the ramification-supported subgroup only;
it does not, by itself, imply that $J(\mathbb Q)[d]=0$.
\end{corollary}

\medskip
\noindent
\textbf{Good reduction and rational torsion.}
The structural results above describe a distinguished part of the
rational $d$-torsion, but they do not control every rational torsion
class. To obtain information about the full group
$J(\mathbb Q)_{\mathrm{tors}}$, we use reduction at primes of good
reduction. The underlying principle is standard: if $p$ is a prime of
good reduction, reduction is injective on the torsion subgroup of
order prime to $p$; see, for example,
Hindry--Silverman \cite[Theorem C.1.4]{HindrySilverman2000}. Thus a rational torsion
point of prime order $\ell\ne p$ forces $\ell$ to divide
$\#J(\mathbb F_p)$. Comparing the orders of reduced Jacobians at
several primes can therefore eliminate possible prime divisors of
rational torsion.

The precise criterion and a convenient sufficient condition are proved
in Section \ref{sec-5}.

\medskip
\noindent
\begin{theorem}
    \label{thm-5.2}[Good-reduction criterion for torsion vanishing]
Let $p_1,\ldots,p_r$ be primes of good reduction for $J$. Suppose that,
for every prime number $\ell$, there is an index $i\in\{1,\ldots,r\}$
such that
\[
p_i\ne\ell
\qquad\text{and}\qquad
\ell\nmid\#J(\mathbb F_{p_i}).
\]
Then
\[
J(\mathbb Q)_{\mathrm{tors}}=0.
\]
\end{theorem}
\medskip
\noindent
\begin{corollary}
    \label{cor-5.3}[A gcd criterion for torsion vanishing]
Let $p_1,\ldots,p_r$ be distinct primes of good reduction for $J$,
where $r\geq2$. Suppose that
\[
\gcd\bigl(\#J(\mathbb F_{p_1}),\ldots,
          \#J(\mathbb F_{p_r})\bigr)=1
\]
and
\[
p_i\nmid\#J(\mathbb F_{p_j})
\qquad\text{whenever }i\ne j.
\]
Then
\[
J(\mathbb Q)_{\mathrm{tors}}=0.
\]
\end{corollary}
The theorem gives a prime-by-prime elimination criterion, while the
corollary packages it into conditions on finitely many group orders.
These results are independent of any assumption that the rational
$d$-torsion is generated by ramification classes. When the orders
$\#J(\mathbb F_p)$ are obtained from point counts or Frobenius
polynomials, the criterion becomes an explicit method for proving
torsion vanishing.

\medskip
\noindent
\textbf{Geometric restrictions from divisors.}
A complementary question is how a divisor class can be represented
on the curve. Fix the rational point $\infty$. Riemann--Roch gives
effective representatives of degree at most $g$ for classes in
$\operatorname{Pic}^0(C)$, while effective divisors of degree $r<g$
map under the Abel--Jacobi map from $\operatorname{Sym}^r(C)$ into a
proper subvariety of $J$. These observations do not classify torsion
points, but they restrict the geometry of classes with small
representatives. The relevant results are established in Section \ref{sec-6}.

\medskip
\noindent
\begin{proposition}
    \label{prop-6.1}[Effective representatives of degree at most
the genus]
Every class $\xi\in\operatorname{Pic}^0(C)$ can be written as
\[
\xi=[E-r\infty],
\]
where $E$ is an effective divisor of degree $r$ for some
$0\leq r\leq g$.
More generally, if \(\xi\) is represented by a divisor \(D\) of degree zero defined over \(\mathbb Q\), then there exist an integer \(r\), with \(0\le r\le g\), and an effective divisor \(E\) of degree \(r\) defined over \(\mathbb Q\) such that
\[
\xi=[E-r\infty].
\]
\end{proposition}

\medskip
\noindent
\begin{proposition}
    \label{prop-6.2}[Finiteness of rational points mapping to
torsion]
Assume $g(C)\geq2$. There are only finitely many points
$P\in C(\mathbb Q)$ for which
\[
[P-\infty]\in J(\mathbb Q)
\]
is a torsion point. This follows from the Manin--Mumford theorem
\cite{Mumford1984}.
\end{proposition}

\medskip
\noindent
\begin{theorem}
    \label{thm-6.3}[Low-degree divisor constraint]
Let $0<r<g$, and let $E_1,E_2$ be distinct effective divisors of
degree $r$. If
\[
[E_1-r\infty]=[E_2-r\infty]\quad\text{in }J,
\]
then the complete linear system $|E_1|$ has positive dimension.
Equivalently,
\[
\ell(E_1)\geq2.
\]
\end{theorem}
The dimension bound for the Abel--Jacobi image and this condition on
linear systems are geometric restrictions rather than a complete
classification of torsion classes. They complement the reduction
criteria by addressing the shape of divisor representatives instead
of the prime divisors of their orders.

\medskip
\noindent
\textbf{Examples and organization.}
Section \ref{sec-7} illustrates how the results interact in explicit examples.

For
\[
C:\quad y^3=x^5-2,
\]
the polynomial \(x^5-2\) is irreducible over \(\mathbb{Q}\), so its
roots form a single Galois orbit. Corollary~\ref{cor-4.4} gives
\(R_3^{G_{\mathbb{Q}}}=0\). This alone does not settle the full
rational torsion, but the finite-field computations at the good
primes \(7\) and \(11\), together with
Corollary~\ref{cor-5.3}, give
\[
J(\mathbb{Q})_{\mathrm{tors}}=0.
\]
For the curve
\[
C:\quad y^3=(x^2-2)(x^2-3),
\]
the roots form two Galois orbits of size \(2\). By
Theorem~\ref{thm-4.2}, the Galois-invariant part of the
ramification subgroup is
\[
R_3^{G_{\mathbb{Q}}}\cong\mathbb{Z}/3\mathbb{Z}.
\]
To determine the full rational \(3\)-torsion in this example, we
use the Prym variety associated with the involution
\((x,y)\mapsto(-x,y)\). The resulting decomposition of the
Jacobian, together with the Galois-module description of Prym
\(3\)-torsion, yields
\[
J(\mathbb{Q})[3]\cong\mathbb{Z}/3\mathbb{Z};
\]
see Section~\ref{subsec:non-transitive-example} and
\cite[\S7.2]{LagaShnidman2025}.

The paper is organized as follows. Section \ref{sec-2} fixes notation and
conventions. Section \ref{sec-3} recalls the background on superelliptic
curves, their ramification, principal divisors, good reduction, and
Weierstrass points. Section \ref{sec-4} proves the structural results for the
ramification subgroup and computes its Galois invariants. Section \ref{sec-5}
establishes the reduction-theoretic criteria for rational torsion
vanishing. Section \ref{sec-6} develops the geometric constraints from effective
divisors and low-degree linear systems. Section \ref{sec-7} gives the two
examples, and Section \ref{sec-8} concludes with further directions.

\subsection*{Acknowledgements}
The author gratefully acknowledges the \textbf{University Grants Commission (UGC), Government of India}, for financial support. The author thanks \textit{\textbf{Prof. Kalyan Banerjee}} for valuable guidance and constructive comments. The author also thanks the \textbf{Department of Mathematics, SRM University-AP}, for providing a supportive research environment.

\subsection*{AI Declaration}

ChatGPT was used solely to assist with the literature review and to check grammatical errors in the manuscript. All mathematical ideas, proofs, calculations, and the writing of the manuscript were carried out independently by the author. The author takes full responsibility for the
mathematical correctness, originality, and content of the paper.

\section{Relevant Notations and Conventions}
\label{sec-2}

Throughout the paper, we use the following notation and conventions.

\begin{itemize}

    \item \textbf{Base field and Galois groups.}
    The base field is $\mathbb{Q}$. We fix an algebraic closure
    $\overline{\mathbb{Q}}$ and write
    \[
        G_{\mathbb{Q}}
        =\operatorname{Gal}
        (\overline{\mathbb{Q}}/\mathbb{Q})
    \]
    for the absolute Galois group. If
    $f(x)\in\mathbb{Z}[x]$ is separable, then
    $\operatorname{Gal}(f/\mathbb{Q})$ denotes the Galois group
    of its splitting field over $\mathbb{Q}$, viewed as a
    permutation group on the roots of $f$. For a prime $p$,
    we write $I_p\subseteq G_{\mathbb{Q}}$ for an inertia
    subgroup at $p$.

    \item \textbf{Group invariants.}
    If a group $G$ acts on an abelian group $M$, we write
    \[
        M^G
        =\{m\in M:\sigma(m)=m
        \text{ for every }\sigma\in G\}
    \]
    for the subgroup of $G$-invariant elements.

    \item \textbf{The superelliptic curve.}
    Let
    \[
        C:\quad y^d=f(x)
    \]
    be a smooth projective curve over $\mathbb{Q}$, where
    $d\geq 3$ is odd and $f(x)\in\mathbb{Z}[x]$ is separable
    of degree $n\geq 3$. We assume $\gcd(d,n)=1$ and write
    $g=g(C)\geq 2$ for the genus. The natural cyclic covering
    is denoted by
    \[
        \pi:C\longrightarrow\mathbb{P}^1,\qquad
        (x,y)\longmapsto x.
    \]
    Under these hypotheses, the covering has degree $d$ and
    is totally ramified above the roots of $f$ and at infinity.

    \item \textbf{The cyclic action.}
    We write $\mu_d$ for the group of $d$-th roots of unity
    in $\overline{\mathbb{Q}}$. The automorphism
    \[
        (x,y)\longmapsto(x,\zeta y),
        \qquad \zeta\in\mu_d,
    \]
    defines an action of $\mu_d$ on
    $C_{\overline{\mathbb{Q}}}$ and hence on its Jacobian.
    This action is understood over $\overline{\mathbb{Q}}$;
    it need not be defined over $\mathbb{Q}$ for every
    $\zeta\in\mu_d$.

    \item \textbf{The Jacobian and its torsion.}
    We denote the Jacobian of $C$ by
    \[
        J=\operatorname{Jac}(C).
    \]
    For a positive integer $m$, the geometric $m$-torsion
    subgroup is
    \[
        J[m]
        =\{P\in J(\overline{\mathbb{Q}}):mP=0\}.
    \]
    The rational $m$-torsion subgroup is
    \[
        J(\mathbb{Q})[m]
        =\{P\in J(\mathbb{Q}):mP=0\}
        =J[m]^{G_{\mathbb{Q}}}.
    \]
    We write $J(\mathbb{Q})_{\mathrm{tors}}$ for the torsion
    subgroup of $J(\mathbb{Q})$. For a prime $\ell$, we denote
    its $\ell$-primary subgroup by
    $J(\mathbb{Q})_{\mathrm{tors}}(\ell)$.

    \item \textbf{Ramification points and the ramification
    subgroup.}
    Let $\alpha_1,\ldots,\alpha_n\in\overline{\mathbb{Q}}$
    be the distinct roots of $f(x)$. For each $i$, let
    \[
        R_i=(\alpha_i,0)\in C(\overline{\mathbb{Q}})
    \]
    be the corresponding ramification point. Since
    $\gcd(d,n)=1$, there is a unique point $\infty$ above
    infinity, and it is defined over $\mathbb{Q}$. We write
    \[
        r_i=[R_i-\infty]\in J(\overline{\mathbb{Q}}).
    \]
    The principal-divisor relation
    \[
        \operatorname{div}(x-\alpha_i)
        =dR_i-d\infty
    \]
    implies that $r_i\in J[d]$. We define the geometric
    ramification subgroup by
    \[
        R_d=\langle r_1,\ldots,r_n\rangle
        \subseteq J[d].
    \]
    The absolute Galois group acts on $R_d$ by permuting the
    ramification classes. In particular,
    \[
        R_d^{G_{\mathbb{Q}}}
        \subseteq J(\mathbb{Q})[d].
    \]
    The inclusion is not assumed to be an equality: the
    ramification subgroup need not account for all rational
    $d$-torsion.

    \item \textbf{Divisors and divisor classes.}
    For a nonzero rational function $h$ on a curve,
    $\operatorname{div}(h)$ denotes its principal divisor.
    We write $D_1\sim D_2$ when two divisors are linearly
    equivalent. Whenever a divisor class is regarded as an
    element of $J$, it is understood to have degree zero.
    If a divisor involves points defined only over
    $\overline{\mathbb{Q}}$, its divisor class is initially
    regarded as a geometric class in
    $J(\overline{\mathbb{Q}})$.

    \item \textbf{Good reduction.}
    A prime $p$ is a prime of good reduction for $C$ if
    $C$ has a smooth proper model over $\mathbb{Z}_p$.
    In this case, $J$ has good reduction at $p$, and its
    special fiber is the Jacobian of the reduced curve over
    $\mathbb{F}_p$. We denote the reduction map on rational
    points by
    \[
        \rho_p:J(\mathbb{Q})\longrightarrow J(\mathbb{F}_p).
    \]
    We use this map to constrain rational torsion by
    comparing the orders of Jacobians over finite fields.

    \item \textbf{Weierstrass points.}
    For a smooth projective curve of genus $g\geq 2$, a
    geometric point $P\in C(\overline{\mathbb{Q}})$ is a
    Weierstrass point if its gap sequence differs from the
    generic sequence $\{1,2,\ldots,g\}$. We write $W(C)$
    for the set of geometric Weierstrass points. This set is
    finite and satisfies the classical bound
    \[
        \#W(C)\leq g^3-g;
    \]
    see \cite{ACGH}.

    \item \textbf{Genus formula.}
    Since $\gcd(d,n)=1$, the Riemann--Hurwitz formula gives
    \[
        g(C)=\frac{(d-1)(n-1)}{2}.
    \]

    \item \textbf{Symmetric powers and the Abel--Jacobi map.}
    For $r\geq 1$, we write $\operatorname{Sym}^r(C)$ for
    the $r$-th symmetric power of $C$, which parametrizes
    effective divisors of degree $r$. Using the rational point
    $\infty$, we define the Abel--Jacobi map
    \[
        \iota:C\longrightarrow J,\qquad
        P\longmapsto[P-\infty].
    \]
    More generally, an effective divisor $E$ of degree $r$
    determines the class $[E-r\infty]\in J$.

\end{itemize}

\section{Background theory}
\label{sec-3}

In this section, we collect the geometric and arithmetic facts used in
the subsequent analysis. We describe the ramification of the cyclic
covering and the resulting genus formula, record the principal divisor
relations associated with the covering, and recall the reduction
property of Jacobians and the classical bound for Weierstrass points.

\subsection{Superelliptic curves and ramification}
\label{subsec:superelliptic-ramification}

Let
\[
    C:\quad y^d=f(x),
\]
where \(f\in\mathbb{Q}[x]\) is separable of degree \(n\), \(d\geq3\)
is odd, and \(\gcd(d,n)=1\). We write \(C\) for the smooth projective
model of this curve and set \(J=\operatorname{Jac}(C)\). The natural
projection
\[
    \pi:C\longrightarrow\mathbb{P}^1,\qquad (x,y)\longmapsto x
\]
is a morphism of degree \(d\).

\begin{proposition}[Ramification structure]
\label{prop-3.1}
Let \(\alpha_1,\ldots,\alpha_n\) be the roots of \(f\) in
\(\overline{\mathbb{Q}}\), and put
\[
    R_i=(\alpha_i,0).
\]
The map \(\pi\) is totally ramified with ramification index \(d\)
above each \(\alpha_i\). Since \(\gcd(d,n)=1\), there is a unique
point \(\infty\in C(\overline{\mathbb{Q}})\) above infinity, and
the map is also totally ramified at that point.
\end{proposition}

\begin{proof}
At a finite point of the affine curve, the equation is
\[
    F(x,y)=y^d-f(x)=0.
\]
Since the curve is smooth in characteristic zero, the ramification
points of the projection to the \(x\)-line occur precisely where
\(\partial F/\partial y=dy^{d-1}=0\). Thus \(y=0\), and the defining
equation then gives \(f(x)=0\). The finite ramification points are
therefore exactly the points \(R_i\).

Because \(f\) is separable, near \(R_i\) we can write
\[
    f(x)=(x-\alpha_i)u(x),\qquad u(\alpha_i)\ne0.
\]
Consequently,
\[
    y^d=(x-\alpha_i)u(x),
\]
so \(y\) is a local parameter at \(R_i\) and
\[
    v_{R_i}(x-\alpha_i)=d.
\]
The ramification index at \(R_i\) is therefore \(d\).

At infinity, the pole orders of \(x\) and \(y\) are determined by
the equation \(y^d=f(x)\). Since \(\gcd(d,n)=1\), there is a single
point above infinity, with pole orders \(d\) for \(x\) and \(n\)
for \(y\). The ramification index there is \(d\). These statements
also follow from the standard ramification theory of superelliptic
curves; see \cite[Chapter III]{Stichtenoth}.
\end{proof}

\begin{proposition}[Genus formula]
\label{prop-3.2}
Under the hypotheses above, the genus of \(C\) is
\[
    g(C)=\frac{(d-1)(n-1)}{2}.
\]
\end{proposition}

\begin{proof}
The morphism \(\pi\) has degree \(d\) and is totally ramified at
the \(n\) finite branch points and at infinity. Applying the
Riemann--Hurwitz formula gives
\[
    2g(C)-2
      =d\bigl(2g(\mathbb{P}^1)-2\bigr)
       +(n+1)(d-1).
\]
Since \(g(\mathbb{P}^1)=0\), this becomes
\[
    2g(C)-2=-2d+(n+1)(d-1),
\]
and hence
\[
    g(C)=\frac{(d-1)(n-1)}{2}.
\]
See also \cite[Chapter III]{Stichtenoth}.
\end{proof}

\begin{lemma}[Good reduction criterion]
\label{lem-3.3}
Suppose that \(f\in\mathbb{Z}[x]\) has degree \(n\), is separable
over \(\mathbb{Q}\), and satisfies \(\gcd(d,n)=1\). Let \(p\) be a
prime such that
\[
    p\nmid d,\qquad
    p\nmid\operatorname{lc}(f),\qquad
    p\nmid\operatorname{disc}(f).
\]
Then the smooth projective model of \(y^d=f(x)\) has good reduction
at \(p\), and the reduced curve has genus
\[
    \frac{(d-1)(n-1)}{2}.
\]
\end{lemma}

\begin{proof}
The discriminant condition implies that the reduction of \(f\)
modulo \(p\) is separable and has the same degree \(n\). Since
\(p\nmid d\), the affine equation
\[
    y^d=f(x)
\]
has no singular points over \(\overline{\mathbb{F}}_p\): a singular
point would satisfy
\[
    y^d=f(x),\qquad f'(x)=0,\qquad dy^{d-1}=0.
\]
The last equation forces \(y=0\), which would give
\(f(x)=f'(x)=0\), contradicting separability.

The leading-coefficient condition preserves the degree of the
equation after reduction. Moreover, \(p\nmid d\) ensures that the
cover is tame, and \(\gcd(d,n)=1\) gives the same single point above
infinity as in characteristic zero. The smooth projective model
therefore has the same ramification structure and genus. These
conditions give the standard good-reduction criterion for this
superelliptic model; see \cite{BouwWewers2017}.
\end{proof}

\subsection{Principal divisors from the covering}
\label{subsec:principal-divisors}

The ramification points give natural torsion classes in the
geometric Jacobian. The relevant relations follow directly from the
divisors of \(x-\alpha_i\) and \(y\).

\begin{proposition}
\label{prop-3.4}
For each root \(\alpha_i\) of \(f\), the following principal-divisor
identity holds on \(C_{\overline{\mathbb{Q}}}\):
\[
    \operatorname{div}(x-\alpha_i)
       =dR_i-d\infty.
\]
Consequently,
\[
    d[R_i-\infty]=0
       \quad\text{in }J(\overline{\mathbb{Q}}).
\]
Moreover,
\[
    \operatorname{div}(y)
       =\sum_{i=1}^{n}R_i-n\infty.
\]
In particular,
\[
    \sum_{i=1}^{n}[R_i-\infty]=0
       \quad\text{in }J(\overline{\mathbb{Q}}).
\]
\end{proposition}

\begin{proof}
On \(\mathbb{P}^1\), the divisor of \(x-\alpha_i\) is
\[
    \operatorname{div}_{\mathbb{P}^1}(x-\alpha_i)
       =(\alpha_i)-(\infty).
\]
Pulling this divisor back under \(\pi\) and using total ramification
at \(R_i\) and \(\infty\) gives
\[
\begin{aligned}
    \operatorname{div}_C(x-\alpha_i)
      &=\pi^*\bigl((\alpha_i)-(\infty)\bigr)\\
      &=dR_i-d\infty.
\end{aligned}
\]
This proves the first identity and the associated \(d\)-torsion
relation.

At each \(R_i\), the function \(y\) vanishes to order one. Its only
pole is at infinity, where it has pole order \(n\). Hence
\[
    \operatorname{div}(y)
       =\sum_{i=1}^{n}R_i-n\infty.
\]
Passing to divisor classes proves the final relation.
\end{proof}

\subsection{Néron models and reduction of torsion}
\label{subsec:reduction-torsion-background}

We recall the reduction property used in Section~\ref{sec-5}.

\begin{theorem}[Injectivity on prime-to-\(p\) torsion]
\label{thm-3.5}
Let \(p\) be a prime of good reduction for \(C\). The reduction map
\[
    \rho_p:J(\mathbb{Q})\longrightarrow J(\mathbb{F}_p)
\]
is injective on the subgroup of rational torsion points whose orders
are relatively prime to \(p\).
\end{theorem}

\begin{proof}
Since \(C\) has good reduction at \(p\), its Jacobian has an abelian
scheme as its good integral model over \(\mathbb{Z}_p\). The kernel
of reduction on \(\mathbb{Q}_p\)-points has no nontrivial torsion
of order prime to \(p\). Thus a torsion point of order prime to
\(p\) that reduces to zero must itself be zero. Restricting this
property to \(J(\mathbb{Q})\) proves the assertion. See
\cite[Theorem C.1.4]{HindrySilverman2000}.
\end{proof}

\subsection{Weierstrass points}
\label{subsec:weierstrass-points}

We finish by recalling the classical bound for Weierstrass points.
This result is background and is not specific to superelliptic curves.

\begin{theorem}[Weierstrass-point bound]
\label{thm-3.6}
Let \(C\) be a smooth projective curve of genus \(g\geq2\) over a
field of characteristic zero. If \(W(C)\) denotes its set of
Weierstrass points over an algebraic closure, then
\[
    \#W(C)\leq g^3-g.
\]
\end{theorem}

\begin{proof}
The bound follows from the classical theory of Weierstrass gap
sequences and the total weight formula for Weierstrass points;
see \cite[Chapter IV]{ACGH}.
\end{proof}

For comparison, a hyperelliptic curve of genus \(g\) has exactly
\(2g+2\) Weierstrass points. The bound above applies more generally
and will only be used as background in the present work.

The divisor identities in
Proposition~\ref{prop-3.4} identify the classes
\([R_i-\infty]\) as natural geometric \(d\)-torsion classes. The
reduction theorem supplies an independent arithmetic constraint on
rational torsion. In the next section, we study the subgroup
generated by the ramification classes and its Galois invariants,
without assuming that this subgroup exhausts the rational
\(d\)-torsion of \(J\).



\section{Ramification Classes and Galois Invariants}
\label{sec-4}

Let
\[
    C:\quad y^d=f(x),\qquad
    f(x)=c\prod_{i=1}^{n}(x-\alpha_i),
\]
where \(f\in\mathbb{Q}[x]\) is separable of degree \(n\),
\(d\geq 3\) is odd, and \(\gcd(d,n)=1\). We work with the
smooth projective model of \(C\) and write \(J=\operatorname{Jac}(C)\).
Let \(R_i=(\alpha_i,0)\) be the finite ramification points, and let
\(\infty\) denote the unique point above infinity. Throughout this
section, the points \(R_i\) and their divisor classes are considered
over \(\overline{\mathbb{Q}}\).

Put
\[
    r_i=[R_i-\infty]\in J[d](\overline{\mathbb{Q}})
\]
and define the ramification subgroup
\[
    R_d=\langle r_1,\ldots,r_n\rangle
    \subseteq J[d](\overline{\mathbb{Q}}).
\]
The purpose of this section is to determine the structure of \(R_d\)
and describe its subgroup of Galois-invariant elements.

\subsection{Principal Divisors and the Ramification Subgroup}
\label{subsec:ramification-structure}

The functions \(x-\alpha_i\) and \(y\) have divisors
\[
    \operatorname{div}(x-\alpha_i)=dR_i-d\infty
    \qquad (1\leq i\leq n)
\]
and
\[
    \operatorname{div}(y)
       =\sum_{i=1}^{n}R_i-n\infty.
\]
Consequently,
\[
    dr_i=0\quad (1\leq i\leq n),
    \qquad
    \sum_{i=1}^{n}r_i=0.
\]

\begin{proof}[Proof of Lemma \ref{lem-4.1}]

Since $\sum_{i=1}^n r_i=0$, we have
$\Delta\subseteq\ker\Phi$. We prove the reverse inclusion.

Let $(a_1,\ldots,a_n)\in\ker\Phi$, and choose representatives
$a_i\in\{0,\ldots,d-1\}$. Then the divisor
\[
A=\sum_{i=1}^n a_iR_i
-\left(\sum_{i=1}^n a_i\right)\infty
\]
is principal. Hence there exists
$h\in\overline{\mathbb{Q}}(C)^\times$ such that
\[
\operatorname{div}(h)=A.
\]

We first observe that the covering
$\pi:C\to\mathbb{P}^1$ is geometrically cyclic of degree
$d$. Indeed, since $f$ is separable, it has a simple root,
say $\alpha_1$. Writing
$f(x)=(x-\alpha_1)g(x)$, with
$g(\alpha_1)\neq 0$, we see that $T^d-f(x)$ is Eisenstein
at $x-\alpha_1$ in $\overline{\mathbb{Q}}[x][T]$.
Consequently,
\[
[\overline{\mathbb{Q}}(C):\overline{\mathbb{Q}}(x)]=d.
\]
Moreover, fixing a primitive $d$-th root of unity
$\zeta_d\in\overline{\mathbb{Q}}$, the automorphism
\[
\sigma(x,y)=(x,\zeta_d y)
\]
generates the Galois group
\[
\operatorname{Gal}\bigl(
\overline{\mathbb{Q}}(C)/\overline{\mathbb{Q}}(x)
\bigr)\cong\mathbb{Z}/d\mathbb{Z}.
\]
In particular, we have the eigenspace decomposition
\[
\overline{\mathbb{Q}}(C)
=\bigoplus_{k=0}^{d-1}
\overline{\mathbb{Q}}(x)y^k,
\]
where $\sigma$ acts on the $k$-th eigenspace by
multiplication by $\zeta_d^k$.

Since $\sigma$ fixes each of the points
$R_1,\ldots,R_n,\infty$, we obtain
\[
\operatorname{div}(\sigma(h))
=\sigma(\operatorname{div}(h))
=A
=\operatorname{div}(h).
\]
Thus $\sigma(h)/h$ has trivial divisor and is therefore
a nonzero constant. Hence $h$ is an eigenfunction for
$\sigma$. By the above eigenspace decomposition, there
exist $u(x)\in\overline{\mathbb{Q}}(x)^\times$ and
$k\in\{0,\ldots,d-1\}$ such that
\[
h=u(x)y^k.
\]

For each $i$, the point $R_i$ is totally ramified under
$\pi$, and hence
\[
v_{R_i}(u(x))\equiv 0\pmod d
\]
for every $u(x)\in\overline{\mathbb{Q}}(x)^\times$.
Since $v_{R_i}(y)=1$, it follows that
\[
a_i=v_{R_i}(h)
=v_{R_i}(u(x))+k
\equiv k\pmod d.
\]
As this holds for every $i$, we conclude that
\[
a_1\equiv a_2\equiv\cdots\equiv a_n\pmod d.
\]
Therefore $(a_1,\ldots,a_n)\in\Delta$, proving that
$\ker\Phi\subseteq\Delta$. Hence
$\ker\Phi=\Delta$, and the first isomorphism theorem
yields
\[
R_d\cong(\mathbb{Z}/d\mathbb{Z})^n/\Delta
\cong(\mathbb{Z}/d\mathbb{Z})^{n-1}.
\]
\end{proof}
\begin{remark}
\label{rem-4.1}
Because \(\gcd(d,n)=1\), every class in
\((\mathbb{Z}/d\mathbb{Z})^n/\Delta\) has a unique representative
\((a_1,\ldots,a_n)\) satisfying
\[
    \sum_{i=1}^{n}a_i=0
       \quad\text{in }\mathbb{Z}/d\mathbb{Z}.
\]
Indeed, adding the diagonal vector \((t,\ldots,t)\) changes the
coordinate sum by \(nt\), and multiplication by \(n\) is invertible
modulo \(d\). In particular, the subgroup \(R_d\) can be represented
uniquely by sum-zero coefficient vectors.
\end{remark}

\subsection{Galois Invariants of the Ramification Subgroup}
\label{subsec:ramification-invariants}

The absolute Galois group
\[
    G_{\mathbb{Q}}
       =\operatorname{Gal}(\overline{\mathbb{Q}}/\mathbb{Q})
\]
permutes the roots of \(f\) and hence the ramification points \(R_i\).
It therefore acts on \(R_d\). Write the decomposition of the roots
into Galois orbits as
\[
    \{\alpha_1,\ldots,\alpha_n\}
       =O_1\sqcup\cdots\sqcup O_s,
    \qquad m_j=|O_j|.
\]

\begin{proof}[Proof of \Cref{thm-4.2}]
By Lemma~\ref{lem-4.1} and
Remark~\ref{rem-4.1}, every element of \(R_d\)
has a unique coefficient-vector representative whose coordinate sum
is zero modulo \(d\). The Galois group permutes the coordinates of
this representative.

Suppose that a class in \(R_d\) is Galois invariant. For every
\(\gamma\in G_{\mathbb{Q}}\), its coefficient vector and the vector
obtained by applying \(\gamma\) represent the same class. Both have
coordinate sum zero, so uniqueness implies that the two vectors are
equal. Hence the coefficients are constant on each Galois orbit
\(O_j\). Write the common coefficient on \(O_j\) as \(b_j\).
The sum-zero condition becomes
\[
    \sum_{j=1}^{s}m_jb_j=0
       \quad\text{in }\mathbb{Z}/d\mathbb{Z}.
\]

Conversely, any vector \((b_1,\ldots,b_s)\) satisfying this
congruence defines a coefficient vector that is constant on each
Galois orbit and whose coordinate sum is zero. It therefore defines
a Galois-invariant class in \(R_d\). The uniqueness of the
sum-zero representative proves that the correspondence is
bijective, and it is plainly a homomorphism.
\end{proof}

\begin{proof}[Proof of Corollary \ref{cor-4.4}]
There is a single Galois orbit, of size \(n\). Theorem
\ref{thm-4.2} therefore reduces the
invariant condition to \(nb=0\) modulo \(d\). Since \(n\) is
invertible modulo \(d\), this forces \(b=0\).
\end{proof}

\subsection{Relation with Rational \texorpdfstring{$d$}{d}-Torsion}
\label{subsec:relation-rational-torsion}

The subgroup \(R_d\) is stable under the Galois action, but its
individual generators need not be rational. Since
\[
    J(\mathbb{Q})[d]
       =J[d](\overline{\mathbb{Q}})^{G_{\mathbb{Q}}},
\]
we have the inclusion
\[
    R_d^{G_{\mathbb{Q}}}
       \subseteq J(\mathbb{Q})[d].
\]
Theorem~\ref{thm-4.2} thus gives an
explicit description of the rational \(d\)-torsion classes that
are generated by the ramification differences.

It does not, by itself, describe all of \(J(\mathbb{Q})[d]\).
A rational \(d\)-torsion class may lie outside \(R_d\). Consequently,
the vanishing of \(R_d^{G_{\mathbb{Q}}}\) does not imply the vanishing
of \(J(\mathbb{Q})[d]\) without an additional argument controlling
the quotient \(J[d](\overline{\mathbb{Q}})/R_d\).

\subsection{The Covering Automorphism and the Base Field}
\label{subsec:covering-automorphism}

The automorphism
\[
    (x,y)\longmapsto(x,\zeta_d y)
\]
is generally defined over \(\mathbb{Q}(\zeta_d)\), rather than over
\(\mathbb{Q}\). Accordingly, a decomposition into character spaces
for the \(\mu_d\)-action must be formulated over a field containing
the \(d\)-th roots of unity, and the action of
\(G_{\mathbb{Q}}\) on those character spaces must be treated
separately.

In particular, individual nontrivial character spaces cannot be
assumed to be \(G_{\mathbb{Q}}\)-stable without proof. No
isotypic-vanishing hypothesis is needed for the results established
in this section: the conclusions above concern only the
ramification-generated subgroup and its Galois invariants.

\section{Reduction-theoretic constraints on rational torsion}
\label{sec-5}

The ramification and Galois-theoretic analysis of the preceding
section describes a subgroup of the rational torsion of the
Jacobian. We now complement that analysis with reduction-theoretic
constraints on the full rational torsion subgroup
\(\Jac(C)(\mathbb Q)_{\mathrm{tors}}\).

The basic principle is that, at a prime of good reduction, the
reduction map is injective on torsion of order prime to the residue
characteristic. Consequently, the order of a rational torsion point
imposes divisibility conditions on the finite groups
\(\Jac(C)(\mathbb F_p)\). Comparing these groups at several primes
can therefore eliminate possible prime divisors of rational torsion.

In Section~\ref{sec-4}, we obtained information about the
Galois-invariant part of the ramification subgroup. The results below
are complementary: they constrain the full rational torsion subgroup
and do not require a hypothesis identifying rational \(d\)-torsion
with the ramification subgroup.

\subsection{Good reduction and torsion constraints}
\label{subsec:good-reduction-torsion}

Let \(J=\Jac(C)\), and let \(p\) be a prime of good reduction for
\(J\). Write
\[
 J(\mathbb Q)_{\mathrm{tors},p'}
 =
 \bigl\{P\in J(\mathbb Q)_{\mathrm{tors}}:
                p\nmid \operatorname{ord}(P)\bigr\}.
\]

The following standard result will be used repeatedly.

\begin{proposition}[Good reduction and prime-to-\(p\) torsion]
\label{prop-5.1}
Let \(p\) be a prime of good reduction for \(J\). The reduction map
\[
 \operatorname{red}_p:
 J(\mathbb Q)_{\mathrm{tors},p'}
 \longrightarrow J(\mathbb F_p)
\]
is injective.
\end{proposition}

\begin{proof}
This is the standard injectivity of reduction on torsion of order
prime to the residue characteristic; see
Hindry--Silverman~\cite[Theorem C.1.4]{HindrySilverman2000}.
\end{proof}

We next record the numerical consequence of this proposition.
Let \(P\in J(\mathbb Q)_{\mathrm{tors}}\) have order
\(N=\operatorname{ord}(P)\). Write
\[
 N=p^{v_p(N)}N_{p'},
 \qquad p\nmid N_{p'},
\]
where \(N_{p'}\) is the prime-to-\(p\) part of \(N\). The point
\[
 P_{p'}=p^{v_p(N)}P
\]
has order \(N_{p'}\). By
Proposition~\ref{prop-5.1}, its reduction has
the same order. It follows that
\[
 N_{p'}\mid \#J(\mathbb F_p).
\]
In particular, if \(\ell\) is a prime divisor of
\(\operatorname{ord}(P)\) and \(\ell\ne p\), then
\[
 \ell\mid \#J(\mathbb F_p).
\]
Thus every prime divisor of the order of a rational torsion point
must divide the order of the reduced Jacobian at every chosen prime
of good reduction different from that torsion prime.

\subsection{Elimination of torsion primes by good reduction}
\label{subsec:elimination-torsion-primes}

The preceding divisibility condition yields a criterion for proving
that the rational torsion subgroup is trivial.

\begin{proof}[Proof of Theorem \ref{thm-5.2}]
Suppose, to the contrary, that
\(J(\mathbb Q)_{\mathrm{tors}}\ne 0\). Choose a nonzero torsion
point \(P\in J(\mathbb Q)_{\mathrm{tors}}\), and let \(\ell\) be
a prime divisor of \(\operatorname{ord}(P)\). The point
\[
 Q=\frac{\operatorname{ord}(P)}{\ell}P
\]
is nonzero and has order \(\ell\).

By hypothesis, there is an index \(i\) such that
\[
 p_i\ne\ell
 \qquad\text{and}\qquad
 \ell\nmid \#J(\mathbb F_{p_i}).
\]
Since \(Q\) has order \(\ell\ne p_i\), its reduction modulo \(p_i\)
has the same order \(\ell\), by
Proposition~\ref{prop-5.1}. Therefore
\[
 \ell\mid \#J(\mathbb F_{p_i}),
\]
a contradiction. Hence
\(J(\mathbb Q)_{\mathrm{tors}}=0\).
\end{proof}

The hypothesis of Theorem~\ref{thm-5.2} can often be
verified by computing the orders of the reduced Jacobians. The
following convenient sufficient condition accounts explicitly for
the possibility that a candidate torsion prime equals one of the
chosen residue characteristics.

\begin{proof}[Proof of Corollary \ref{cor-5.3}]
We verify the hypothesis of
Theorem~\ref{thm-5.2} for every prime \(\ell\).

First suppose that
\(\ell\notin\{p_1,\ldots,p_r\}\). If \(\ell\) divided
\(\#J(\mathbb F_{p_i})\) for every \(i\), then it would divide their
gcd, contradicting the assumption that this gcd is \(1\). Hence
there is an index \(i\) such that
\[
 \ell\nmid \#J(\mathbb F_{p_i}).
\]
Since \(\ell\ne p_i\), this is the required index.

Now suppose that \(\ell=p_k\) for some \(k\). Since \(r\ge2\),
choose \(j\ne k\). The second hypothesis gives
\[
 \ell=p_k\nmid \#J(\mathbb F_{p_j}),
\]
and \(p_j\ne\ell\). Thus the required index exists in this case
as well. Theorem~\ref{thm-5.2} now gives the result.
\end{proof}

\subsection{Application to the superelliptic family}
\label{subsec:reduction-superelliptic-family}

Consider
\[
 C:\quad y^d=f(x),
\]
where \(d\ge3\) is odd, \(f(x)\in\mathbb Q[x]\) is separable of
degree \(n\ge3\), and
\[
 \gcd(d,n)=1.
\]
Let \(J=\Jac(C)\). Under these hypotheses, the smooth projective
curve has a unique point at infinity, and its genus is
\[
 g=\frac{(d-1)(n-1)}{2}.
\]

To apply the reduction criteria, choose an integral model of the
curve and a set \(\mathcal P_{\mathrm{good}}\) of primes at which
that model has good reduction. The set of excluded primes depends
on the chosen integral model. In particular, for a model
\(y^d=f(x)\) with \(f(x)\in\mathbb Z[x]\), primes dividing \(d\),
the leading coefficient of \(f\), or the discriminant of \(f\)
must be excluded when using the usual smooth-model criterion.
More generally, denominators introduced when passing from a
rational equation to an integral model must also be taken into
account. The good-reduction hypothesis should be verified for
the model actually used.

For \(p\in\mathcal P_{\mathrm{good}}\), let \(C_p\) and \(J_p\)
denote the reductions of \(C\) and \(J\), respectively. Set
\[
 N_p=\#J_p(\mathbb F_p)=\#J(\mathbb F_p).
\]
These finite-field cardinalities are the numerical input for
Theorem~\ref{thm-5.2} and
Corollary~\ref{cor-4.4}.

Indeed, suppose that
\(P\in J(\mathbb Q)_{\mathrm{tors}}\) is nonzero, and let
\(\ell\) be a prime divisor of \(\operatorname{ord}(P)\). Then
\[
 Q=\frac{\operatorname{ord}(P)}{\ell}P
\]
has order \(\ell\). For any
\(p\in\mathcal P_{\mathrm{good}}\) with \(p\ne\ell\), the
reduction of \(Q\) has order \(\ell\), so
\[
 \ell\mid N_p.
\]
Consequently, to prove that the rational torsion subgroup is
trivial, it is enough to choose good primes
\(p_1,\ldots,p_r\) and show that for every prime \(\ell\) there
is an index \(i\) such that
\[
 p_i\ne\ell
 \qquad\text{and}\qquad
 \ell\nmid N_{p_i}.
\]
In practice, one can factor the integers \(N_{p_i}\) and use these
factorizations to eliminate candidate torsion primes.

The cardinalities \(N_p\) can be computed from Frobenius data.
With the standard normalization
\[
 P_p(T)=\prod_{j=1}^{2g}(1-\alpha_jT),
\]
where the \(\alpha_j\) are the Frobenius eigenvalues on the first
étale cohomology of \(C_p\), the characteristic polynomial
\(P_p(T)\) satisfies
\[
 \#J_p(\mathbb F_p)=P_p(1).
\]
Equivalently, the required cardinality can be obtained from the
Frobenius polynomial of the curve or Jacobian over \(\mathbb F_p\).
The precise convention for arithmetic versus geometric Frobenius
must be kept consistent with the polynomial normalization.

The reduction-theoretic information developed here complements
the ramification analysis in
Section~\ref{sec-4}. That section computes the
Galois-invariant part of the subgroup generated by the ramification
classes and places constraints on rational \(d\)-torsion arising
from the geometry of the covering. The present section instead
constrains the full rational torsion subgroup by using reductions
at suitable good primes. It does not rely on an equality between
the ramification subgroup and all of \(J[d]\).

The two approaches can therefore be used together. Structural
information about \(J(\mathbb Q)[d]\), when available, can be
combined with the reduction constraints to eliminate possible
torsion primes. If the resulting arguments rule out every prime
divisor of the order of a rational torsion point, then
\[
 J(\mathbb Q)_{\mathrm{tors}}=0.
\]
The examples below illustrate how finite-field computations can
supply the reduction data needed for these criteria.

\section{Geometric constraints on exceptional torsion}
\label{sec-6}

The reduction-theoretic results of
Section~\ref{sec-5} constrain the possible prime
divisors of \(J(\mathbb Q)_{\mathrm{tors}}\) through good reduction.
We now consider complementary geometric restrictions arising from
the degrees of divisors representing classes in the Jacobian.

The key observation is that classes admitting effective
representatives of small degree are controlled by the symmetric
powers of \(C\) and their images under the Abel--Jacobi map. Since
the \(r\)-th symmetric power has dimension \(r\), these images have
dimension at most \(r\). In particular, when \(r<g=\dim J\), they
lie in proper subvarieties of \(J\). We also explain how equality
between classes represented by distinct effective divisors imposes
a condition on the associated complete linear system.

Throughout this section, \(C\) is a smooth projective curve over
\(\mathbb Q\) of genus \(g\), \(J=\operatorname{Jac}(C)\), and
\(\infty\in C(\mathbb Q)\) is a fixed rational point. In our
superelliptic setting, the coprimality condition
\(\gcd(d,\deg f)=1\) supplies such a unique point at infinity.

\subsection{Reduced representatives and Abel--Jacobi constraints}
\label{subsec:reduced-representatives}

We first recall the existence of effective representatives of
degree at most \(g\).

\begin{proof}[Proof of Proposition \ref{prop-6.1}]

Let \(D\) be a divisor of degree zero over \(\mathbb Q\) representing \(\xi\). By the Riemann–Roch theorem,
\[
\ell(D+g\infty)-\ell(K-D-g\infty)
=\deg(D+g\infty)-g+1=1,
\]
where \(K\) is a canonical divisor. Since \(\ell(K-D-g\infty)\ge0\), we obtain
\[
\ell(D+g\infty)\ge1.
\]
Thus there exists a nonzero rational function \(h\in\mathbb Q(C)\) such that
\[
E=D+g\infty+\operatorname{div}(h)
\]
is effective. The divisor \(E\) is defined over \(\mathbb Q\) and has degree \(g\). Therefore,
\[
[D]=[E-g\infty],
\]
which proves the first assertion.

For the more general statement, let \(r\) be the smallest nonnegative integer such that
\[
\ell(D+r\infty)>0.
\]
Such an integer exists by the first part, and \(r\le g\). Choose a nonzero \(h\in\mathbb Q(C)\) such that
\[
E=D+r\infty+\operatorname{div}(h)
\]
is effective. By the minimality of \(r\), the effective divisor \(E\) cannot contain \(\infty\); otherwise, subtracting \(\infty\) would give a nonzero section of \(D+(r-1)\infty\) when \(r>0\). Hence \(E\) has degree \(r\), and
\[
[D]=[E-r\infty].
\]
This completes the proof. 

\end{proof}
The proposition gives an effective representative of degree \(g\).
It does not assert that every class has a representative of degree
strictly less than \(g\). Whenever a class has a representative
\([E-r\infty]\) with \(r<g\), the Abel--Jacobi map places that class
in a proper subvariety of \(J\), as explained below.

\begin{proof}[Proof of Proposition \ref{prop-6.2}]
    
 Consider the Abel–Jacobi embedding
\[
\iota:C\longrightarrow J,\qquad
P\longmapsto[P-\infty].
\]
Since \(g\ge2\), the image \(\iota(C)\) is not a translate of a positive-dimensional abelian subvariety of \(J\). By the Manin–Mumford theorem \cite{Mumford1984}, the intersection of \(\iota(C)\) with the set \(J(\overline{\mathbb Q})_{\mathrm{tors}}\) of geometric torsion points is finite. Consequently, its subset consisting of images of rational points is also finite. This proves the proposition.

\end{proof}

\subsection{Low-degree linear systems}
\label{subsec:low-degree-linear-systems}

The next result gives a necessary condition for two distinct
effective divisors of small degree to determine the same class in
the Jacobian.

\begin{proof}[Proof of \Cref{thm-6.3}]

The equality of divisor classes gives
\[
E_1-r\infty\sim E_2-r\infty,
\]
and hence \(E_1\sim E_2\). Thus \(E_1\) and \(E_2\) are distinct effective members of the complete linear system \(|E_1|\). Therefore \(|E_1|\) contains at least two distinct members, so
\[
\dim |E_1|\ge1.
\]
Since \(\dim |E_1|=\ell(E_1)-1\), we conclude that \(\ell(E_1)\ge2\).
\end{proof}

\begin{remark}
\label{rem-6.1}
For \(0<r<g\), effective divisors of degree \(r\) are parametrized by the
symmetric power \(\operatorname{Sym}^{r}(C)\), which has dimension \(r\).
The Abel--Jacobi map
\[
u_r:\operatorname{Sym}^{r}(C_{\overline{\mathbb{Q}}})
\longrightarrow J_{\overline{\mathbb{Q}}},
\qquad
E\longmapsto [E-r\infty]
\]
therefore has image of dimension at most \(r<g=\dim J\). Consequently,
the geometric classes admitting effective representatives of degree
\(r<g\) lie in a proper subvariety of \(J_{\overline{\mathbb{Q}}}\).

This dimension bound is a geometric restriction; it does not imply that
the image contains only finitely many torsion points. Theorem \ref{thm-6.3}
provides an additional restriction: if two distinct effective divisors
of the same degree \(r<g\) represent the same divisor class, then their
common complete linear system has positive dimension. These observations
constrain the geometry of divisor representatives but do not, by
themselves, classify the rational torsion subgroup \(J(\mathbb{Q})_{\mathrm{tors}}\).
\end{remark}

The geometric constraints in this section complement the
reduction-theoretic criteria of
Section~\ref{sec-5}. Good reduction restricts
possible prime divisors of rational torsion through the finite
groups \(J(\mathbb F_p)\), whereas the Abel--Jacobi map and the
geometry of linear systems constrain classes admitting
representatives of small degree. Neither method alone establishes
a general classification of \(J(\mathbb Q)_{\mathrm{tors}}\), but
they provide distinct restrictions that may be combined when
studying a specific curve.

\section{Examples}
\label{sec-7}

We illustrate the structural, reduction-theoretic, and geometric
constraints developed in the preceding sections through two explicit
superelliptic curves. The first example has a transitive Galois action
on the roots of the defining polynomial and admits a reduction-theoretic
argument for the vanishing of the full rational torsion subgroup.
The second example has two Galois orbits of roots and illustrates how
the orbit decomposition determines the Galois-invariant part of the
ramification subgroup.

\subsection{A transitive genus-\texorpdfstring{$4$}{4} example}
\label{subsec:transitive-genus-four}

Consider the superelliptic curve
\[
  C:\quad y^3=x^5-2
\]
over \(\mathbb Q\), and let \(J=\operatorname{Jac}(C)\).

Since \(\gcd(3,5)=1\), the genus formula gives
\[
  g(C)=\frac{(3-1)(5-1)}{2}=4.
\]
Moreover, \(f(x)=x^5-2\) is Eisenstein at \(2\), and hence
irreducible over \(\mathbb Q\). Thus the Galois group acts
transitively on its five roots. There is a single Galois orbit
of size \(m_1=5\).

By Theorem~\ref{thm-4.2}, the
Galois-invariant part of the ramification subgroup is described by
\[
  R_3^{G_{\mathbb Q}}
  \cong
  \{a\in\mathbb Z/3\mathbb Z:5a\equiv0\pmod3\}.
\]
Since \(5\equiv2\pmod3\), the only solution is \(a=0\). Hence
\[
  R_3^{G_{\mathbb Q}}=0.
\]
This conclusion concerns the ramification subgroup only; it does
not, by itself, imply that \(J(\mathbb Q)[3]=0\).

We next use good reduction to constrain the full rational torsion
subgroup. The discriminant of \(f(x)=x^5-2\) is
\[
  \operatorname{disc}(f)
  =(-1)^{5(5-1)/2}5^5(-2)^4
  =2^4 5^5.
\]
The leading coefficient of \(f\) is \(1\). Therefore, the primes
\(p=7\) and \(p=11\) do not divide
\(3\operatorname{lc}(f)\operatorname{disc}(f)\). Both are primes
of good reduction for the curve.

For \(p\in\{7,11\}\) and \(k\ge1\), put
\[
  N_k(p)=\#C(\mathbb F_{p^k}).
\]
Since \(\gcd(3,5)=1\), the smooth projective model has a unique
point at infinity. Thus the projective point count is obtained by
counting the affine solutions of
\[
  y^3=x^5-2
\]
over \(\mathbb F_{p^k}\) and adding this point at infinity.

The point counts reported in the manuscript are
\[
\begin{array}{c|rrrr}
 p & N_1(p)&N_2(p)&N_3(p)&N_4(p)\\
 \hline
 7  & 8  & 50  & 344  & 2766\\
 11 & 12 & 154 & 1332 & 15274
\end{array}
\]
These counts are stated to have been obtained by exhaustive
finite-field enumeration and checked by the same procedure. The
SageMath verification is included in Appendix~\ref{app}.

We now recover the Frobenius polynomials from these counts. For a
prime \(p\) of good reduction, let
\[
  P_p(T)=\prod_{i=1}^{8}(1-\alpha_iT),
\]
where the \(\alpha_i\) are the Frobenius eigenvalues, with the
normalization for which
\[
  \#J(\mathbb F_p)=P_p(1).
\]
The point-counting formula gives
\[
  N_k(p)=p^k+1-\sum_{i=1}^{8}\alpha_i^k.
\]
Writing
\[
  S_k(p)=p^k+1-N_k(p),
\]
the values \(S_k(p)\) are the power sums of the Frobenius
eigenvalues. Newton's identities determine the first four
coefficients of \(P_p(T)\), and the functional equation determines
the remaining coefficients.

For \(p=7\), the point counts give
\[
  S_1(7)=0,\qquad S_2(7)=0,\qquad
  S_3(7)=0,\qquad S_4(7)=-364.
\]
Newton's identities give
\[
  e_1=e_2=e_3=0,\qquad e_4=91,
\]
where \(e_i\) denotes the \(i\)-th elementary symmetric function
of the Frobenius eigenvalues. The functional equation then yields
\[
  P_7(T)=1+91T^4+2401T^8.
\]
Consequently,
\[
  \#J(\mathbb F_7)=P_7(1)=1+91+2401=2493.
\]

Similarly, for \(p=11\),
\[
  S_1(11)=0,\qquad S_2(11)=-32,\qquad
  S_3(11)=0,\qquad S_4(11)=-632.
\]
Newton's identities give
\[
  e_1=0,\qquad e_2=16,\qquad
  e_3=0,\qquad e_4=286.
\]
Using the functional equation, we obtain
\[
  P_{11}(T)
  =1+16T^2+286T^4+1936T^6+14641T^8.
\]
It follows that
\[
\begin{aligned}
  \#J(\mathbb F_{11})
  &=P_{11}(1)\\
  &=1+16+286+1936+14641\\
  &=16880.
\end{aligned}
\]

We now verify the hypotheses of
Corollary~\ref{cor-5.3}. First,
\[
  \gcd\bigl(\#J(\mathbb F_7),\#J(\mathbb F_{11})\bigr)
  =\gcd(2493,16880)=1.
\]
Moreover,
\[
  7\nmid16880
  \qquad\text{and}\qquad
  11\nmid2493.
\]
Thus the primes \(7\) and \(11\) satisfy both the gcd and the
cross-divisibility conditions of the corollary. Consequently,
\[
  \boxed{J(\mathbb Q)_{\mathrm{tors}}=0.}
\]
This conclusion follows from the reduction-theoretic calculation
and does not depend on any additional hypothesis concerning the
ramification subgroup.

\subsection{A non-transitive example}
\label{subsec:non-transitive-example}

Consider the curve
\[
    C:\quad y^3=(x^2-2)(x^2-3)
       =x^4-5x^2+6
\]
over \(\mathbb{Q}\), and let \(J\) denote its Jacobian. The
polynomial
\[
    f(x)=(x^2-2)(x^2-3)
\]
has two Galois orbits of roots, namely
\[
    \{\sqrt{2},-\sqrt{2}\}
    \qquad\text{and}\qquad
    \{\sqrt{3},-\sqrt{3}\}.
\]
Consequently, Theorem~\ref{thm-4.2} gives
\[
    R_3^{G_{\mathbb{Q}}}
    \cong
    \left\{(a,b)\in(\mathbb{Z}/3\mathbb{Z})^2:
                 2a+2b=0\right\}
    \cong \mathbb{Z}/3\mathbb{Z}.
\]
Thus the ramification subgroup supplies a nontrivial rational
subgroup of order \(3\). We now determine the full rational
\(3\)-torsion of \(J\).

Consider the involution
\[
    \iota:C\longrightarrow C,\qquad
    (x,y)\longmapsto(-x,y).
\]
Its quotient is the genus-one curve
\[
    D:\quad y^3=u^2-5u+6,\qquad u=x^2.
\]
The curve \(D\) has a rational point, so its Jacobian is an
elliptic curve over \(\mathbb{Q}\). The change of variables
\[
    X=4y,\qquad Y=8u-20
\]
transforms the equation of \(D\) into the Weierstrass equation
\[
    E:\quad Y^2=X^3+16.
\]
Consequently, the Jacobian of \(D\) is isomorphic over
\(\mathbb{Q}\) to \(E\).

The point \((0,4)\in E(\mathbb{Q})\) has order \(3\). Moreover,
the \(3\)-division polynomial of \(E\) is
\[
    \psi_3(X)=3X(X^3+64).
\]
Its only rational roots are \(X=0\) and \(X=-4\). The latter
does not give a rational point on \(E\), since it yields
\(Y^2=-48\). For \(X=0\), we obtain precisely the points
\((0,\pm4)\). It follows that
\[
    E(\mathbb{Q})[3]\cong\mathbb{Z}/3\mathbb{Z}.
\]

We next study the Prym variety associated with \(\iota\). Write
\(P\) for the Prym variety of the double cover \(C\to D\), so
that
\[
    P=\ker(1+\iota^*:J\longrightarrow J)^0.
\]
The curve is a member of the bielliptic Picard family
\[
    C_{a,b}:\quad y^3=x^4+ax^2+b
\]
with \(a=-5\) and \(b=6\). The smoothness condition in this
family is satisfied, since
\[
    16b(a^2-4b)=16\cdot6\cdot(25-24)=96\ne0.
\]
Thus the Prym results of Laga--Shnidman
\cite[\S2.5, \S7.2]{LagaShnidman2025} apply to this curve.
In particular, the natural map
\[
    E\times P\longrightarrow J
\]
is an isogeny with kernel of order a power of \(2\). It therefore
induces an isomorphism on rational \(3\)-primary torsion, and
\[
    J(\mathbb{Q})[3]\cong
    E(\mathbb{Q})[3]\oplus P(\mathbb{Q})[3].
\]
It remains to show that \(P(\mathbb{Q})[3]=0\).

Let \(\omega\) be a primitive cube root of unity and put
\(\mathfrak p=(1-\omega)\) in \(\mathbb{Z}[\omega]\). By
Lemma~7.8 of \cite{LagaShnidman2025}, the group
\(P[\mathfrak p]\) is a two-dimensional vector space over
\(\mathbb{F}_3\), with basis
\[
    D_1=(\sqrt{2},0)-(-\sqrt{2},0),\qquad
    D_2=(\sqrt{3},0)-(-\sqrt{3},0).
\]
Here the points are viewed as points on the smooth projective
model of \(C\), and the differences denote their corresponding
divisor classes in the Prym.

By Proposition~7.9 of \cite{LagaShnidman2025}, the Galois action
on this basis is induced by the action on the root pairs. More
explicitly, changing the sign of \(\sqrt{2}\) while fixing
\(\sqrt{3}\) sends \(D_1\) to \(-D_1\) and fixes \(D_2\);
changing the sign of \(\sqrt{3}\) while fixing \(\sqrt{2}\)
fixes \(D_1\) and sends \(D_2\) to \(-D_2\). Consequently,
the representation is
\[
    \rho(\sigma)=
    \begin{pmatrix}
        \chi_2(\sigma)&0\\
        0&\chi_3(\sigma)
    \end{pmatrix},
\]
where \(\chi_2\) and \(\chi_3\) are the quadratic characters
associated with \(\mathbb{Q}(\sqrt{2})\) and
\(\mathbb{Q}(\sqrt{3})\), respectively. Since both characters
are nontrivial, there are no nonzero Galois-invariant vectors,
and hence
\[
    P[\mathfrak p](\mathbb{Q})=0.
\]

Equation~(7.3) of \cite{LagaShnidman2025} gives an exact sequence
of Galois modules
\[
    0\longrightarrow P[\mathfrak p]
    \longrightarrow P[3]
    \longrightarrow P_{-27}[\mathfrak p]
    \longrightarrow0,
\]
where \(P_{-27}\) is the quadratic twist of \(P\) by
\(\mathbb{Q}(\sqrt{-3})/\mathbb{Q}\). Under this twist, the two
characters describing the action on the corresponding
\(\mathfrak p\)-torsion become
\[
    \chi_2\chi_{-3}=\chi_{-6},
    \qquad
    \chi_3\chi_{-3}=\chi_{-9}=\chi_{-1}.
\]
Both characters are nontrivial. Therefore
\[
    P_{-27}[\mathfrak p](\mathbb{Q})=0.
\]
Taking Galois invariants in the exact sequence, we conclude that
\[
    P[3](\mathbb{Q})=0,
\]
or equivalently \(P(\mathbb{Q})[3]=0\). It follows that
\[
    \boxed{J(\mathbb{Q})[3]\cong\mathbb{Z}/3\mathbb{Z}.}
\]
In this example, the rational \(3\)-torsion is therefore exactly
the subgroup supplied by the Galois-invariant ramification classes.

\section{Concluding remarks and further directions}
\label{sec-8}

In this paper, we have studied rational torsion on the Jacobians of odd
superelliptic curves
\[
    C:\quad y^d=f(x),
\]
where \(f\in\mathbb{Q}[x]\) is separable, \(d\geq 3\) is odd, and
\(\gcd(d,\deg f)=1\). Our approach combines the geometry of the
ramification points of the covering \(C\to\mathbb{P}^1\), the Galois
action on these points, and reduction-theoretic constraints on rational
torsion.

The structural results of Section~\ref{sec-4} describe the
subgroup \(R_d\subseteq J[d](\overline{\mathbb{Q}})\) generated by the
classes of the ramification points relative to the point at infinity.
They also give an explicit description of its Galois-invariant
subgroup in terms of the Galois orbits of the roots of \(f\). In
particular,
\[
    R_d^{G_{\mathbb{Q}}}\subseteq J(\mathbb{Q})[d],
\]
and the orbit structure can therefore yield concrete information about
rational \(d\)-torsion. This calculation concerns the
ramification-generated subgroup; without additional arguments, it
does not determine all of \(J(\mathbb{Q})[d]\).

The reduction-theoretic results of Section~\ref{sec-5}
provide a complementary method. At primes of good reduction, the
prime-to-\(p\) torsion of \(J(\mathbb{Q})\) injects into
\(J(\mathbb{F}_p)\). Consequently, suitable finite-field point counts
can eliminate individual prime-primary components of rational torsion.
Under the hypotheses of
Corollary~\ref{cor-5.3}, these constraints imply
\[
    J(\mathbb{Q})_{\mathrm{tors}}=0.
\]
Thus, the reduction method can establish complete torsion vanishing
without requiring a separate hypothesis about the ramification-generated
subgroup.

The examples in Section~\ref{sec-7} illustrate how the
Galois-invariant ramification subgroup and reduction-theoretic
methods complement one another. For the curve
\[
    C:\quad y^3=x^5-2,
\]
the transitive Galois action on the five roots implies
\[
    R_3^{G_{\mathbb{Q}}}=0.
\]
The point counts at the good primes \(7\) and \(11\), together
with the required cross-divisibility conditions, yield
\[
    J(\mathbb{Q})_{\mathrm{tors}}=0
\]
by Corollary~\ref{cor-5.3}.

For the curve
\[
    C:\quad y^3=(x^2-2)(x^2-3),
\]
the two Galois orbits of roots give
\[
    R_3^{G_{\mathbb{Q}}}\cong\mathbb{Z}/3\mathbb{Z}.
\]
Thus the ramification classes provide a nontrivial rational
subgroup of order \(3\). The Prym variety associated with the
involution \((x,y)\mapsto(-x,y)\) allows us to go further:
the quotient curve has Jacobian isomorphic to
\[
    E:\quad Y^2=X^3+16,
\]
and the Prym \(3\)-torsion calculation shows that the Prym
contributes no nontrivial rational \(3\)-torsion. Consequently,
\[
    J(\mathbb{Q})[3]\cong\mathbb{Z}/3\mathbb{Z}.
\]
This example therefore illustrates how combining the orbit
description of ramification classes with the arithmetic of a
Prym variety can determine the full rational \(3\)-torsion,
rather than only the contribution from the ramification subgroup;
see \cite[\S7.2]{LagaShnidman2025}.

Several questions merit further investigation.

\begin{enumerate}
    \item[\textup{(1)}] \textbf{Torsion outside the ramification subgroup.}
    The quotient \(J[d]/R_d\) records the part of the geometric
    \(d\)-torsion not accounted for by the ramification classes.
    Understanding its Galois action, together with the cohomological
    obstruction associated with
    \[
        0\longrightarrow R_d\longrightarrow J[d]
        \longrightarrow J[d]/R_d\longrightarrow 0,
    \]
    may yield additional criteria for determining
    \(J(\mathbb{Q})[d]\). In particular, one should identify hypotheses
    that genuinely imply the desired description of rational torsion,
    rather than infer it from the vanishing of quotient invariants alone.

    \item[\textup{(2)}] \textbf{Non-coprime degrees.}
    It would be useful to extend the structural analysis to\\
    \(\gcd(d,\deg f)\ne 1\). In this setting, the points above infinity
    and their divisor classes require separate treatment, and the
    ramification relations may differ from those considered here.

    \item[\textup{(3)}] \textbf{Curves over general number fields.}
    The structural and reduction-theoretic arguments can be investigated
    over a general number field \(K\). Such an extension would require
    accounting for the \(G_K\)-action, the residue fields at primes of
    good reduction, and the corresponding reduction maps.

    \item[\textup{(4)}] \textbf{Finite-field computations.}
    Systematic computation of Frobenius polynomials and Jacobian orders
    for explicit superelliptic families may provide effective ways to
    select good primes that eliminate rational torsion. The interaction
    between these computations and the automorphisms of the covering
    may also be worth exploring.

    \item[\textup{(5)}] \textbf{Low-degree divisor geometry.}
    The results of Section~\ref{sec-6} give necessary
    geometric conditions for distinct effective divisors of degree
    less than the genus to represent the same divisor class. A more
    detailed study of the Abel--Jacobi images of
    \(\operatorname{Sym}^r(C)\), for \(r<g\), may provide further
    restrictions on rational torsion classes admitting effective
    representatives of small degree.
\end{enumerate}

Overall, the results show that ramification geometry, Galois symmetry,
and arithmetic reduction provide complementary tools for studying
rational torsion on odd superelliptic Jacobians. Determining when these
tools can be combined into general classification theorems remains an
interesting problem.

\appendix

\section{Finite-field verification for Example~\ref{subsec:transitive-genus-four}}
\label{app}

We verify the finite-field point counts for the curve
\[
    C:\quad y^3=x^5-2
\]
over \(\mathbb{F}_{p^k}\), where \(p\in\{7,11\}\) and
\(1\leq k\leq4\). Since \(\gcd(3,5)=1\), the smooth projective
model has a unique point at infinity. Thus, if \(N_k\) denotes the
number of \(\mathbb{F}_{p^k}\)-rational points on the projective
curve, then
\[
    N_k=1+
    \#\bigl\{(x,y)\in\mathbb{F}_{p^k}^2:
                     y^3=x^5-2\bigr\}.
\]

The following SageMath code counts the cube roots of each field
element once and then uses these counts to enumerate the affine
points. The assertions check the values used in
Subsection~\ref{subsec:transitive-genus-four}.

\begin{verbatim}
def point_count(p, k):
    F = GF(p^k, name='a')

    # Count the cube roots of each element of F.
    cube_counts = {}
    for y in F:
        z = y^3
        cube_counts[z] = cube_counts.get(z, 0) + 1

    # Count affine points satisfying y^3 = x^5 - 2.
    count_affine = 0
    for x in F:
        rhs = x^5 - F(2)
        count_affine += cube_counts.get(rhs, 0)

    # There is exactly one point at infinity.
    return count_affine + 1


expected = {
    7:  [8, 50, 344, 2766],
    11: [12, 154, 1332, 15274],
}

for p in [7, 11]:
    computed = [point_count(p, k) for k in [1, 2, 3, 4]]
    print("p =", p, "point counts =", computed)
    assert computed == expected[p]
\end{verbatim}

The expected output is
\begin{verbatim}
p = 7 point counts = [8, 50, 344, 2766]
p = 11 point counts = [12, 154, 1332, 15274]
\end{verbatim}

\printbibliography

\end{document}